\documentclass[12pt, reqno]{amsart}
\usepackage[dvips]{graphicx}
\usepackage[arrow,matrix,curve]{xy}

\newtheorem{theorem}{Theorem}[section]
\newtheorem{proposition}[theorem]{Proposition}
\newtheorem{lemma}[theorem]{Lemma}
\newtheorem{corollary}[theorem]{Corollary}
\newtheorem{conjecture}[theorem]{Conjecture}

\theoremstyle{definition} 
\newtheorem{definition}[theorem]{Definition}
\newtheorem{example}[theorem]{Example}
\newtheorem{remark}[theorem]{Remark}
\newtheorem{question}[theorem]{Question}

\newcommand{\C}{\mathbb C} 
\newcommand{\R}{\mathbb R} \newcommand{\T}{\mathbb T}
\newcommand{\Z}{\mathbb Z}

\DeclareMathOperator{\dist}{{\rm dist}}
\DeclareMathOperator{\sys}{{\rm sys}}
\DeclareMathOperator{\stsys}{{\rm stsys}}
\DeclareMathOperator{\confsys}{{\rm confsys}}
\DeclareMathOperator{\pisys}{{\rm sys}\pi}

\DeclareMathOperator{\vol}{{\rm vol}}
\DeclareMathOperator{\hyperelliptic}{{\it J}}
\DeclareMathOperator{\tauzero}{{\tau_0}}
\DeclareMathOperator{\FillVol}{{\rm FillVol}}

\DeclareMathOperator{\area}{{\rm area}} 

\def\length{{\rm length}}

\def\AF {{AF(\pi,\pi)}}

\def\gmetric {{\mathcal G}}

\def\genus {{\it g}}

\def\ie {{\it i.e.\ }} 
 
\def\cf {\hbox{\it cf.\ }}

\def\Ci{{C_0}}

\numberwithin{equation}{section} 
\numberwithin{figure}{section}

\begin{document}

\author[V.~Bangert]{Victor Bangert$^!$} \address{ Mathematisches
Institut, Universit\"at Freiburg, Eckerstr.~1, 79104 Freiburg,
Germany} \email{bangert@mathematik.uni-freiburg.de}
\thanks{$^!$Partially Supported by DFG-Forschergruppe `Nonlinear
Partial Differential Equations: Theoretical and Numerical Analysis'}

\author[C.~Croke]{Christopher Croke$^+$} \address{ Department of
Mathematics, University of Pennsylvania, Philadelphia, PA 19104-6395
USA} \email{ccroke@math.upenn.edu} \thanks{$^+$Supported by NSF grant
DMS 02-02536 and the MSRI}

\author[S.~Ivanov]{Sergei V. Ivanov$^{\dagger}$} \address{Steklov
Math. Institute, Fontanka 27, RU-191011 St. Petersburg, Russia}
\email{svivanov@pdmi.ras.ru} \thanks{$^\dagger$Supported by grants
CRDF RM1-2381-ST-02, RFBR 02-01-00090, and NS-1914.2003.1}

\author[M.~Katz]{Mikhail G. Katz$^{*}$} \address{Department of
Mathematics, Bar Ilan University, Ramat Gan 52900
Israel} \email{katzmik@math.biu.ac.il} \thanks{$^{*}$Supported by the
Israel Science Foundation (grants no.\ 620/00-10.0 and 84/03)}

\title[Filling area conjecture and hyperelliptic surfaces] {Filling
area conjecture and ovalless real hyperelliptic surfaces$^{1}$}



\begin{abstract}
We prove the filling area conjecture in the hyperelliptic case.  In
particular, we establish the conjecture for all genus~1 fillings of
the circle, extending P. Pu's result in genus 0.  We translate the
problem into a question about closed ovalless real surfaces.  The
conjecture then results from a combination of two ingredients.  On the
one hand, we exploit integral geometric comparison with orbifold
metrics of constant positive curvature on real surfaces of even
positive genus.  Here the singular points are Weierstrass points.  On
the other hand, we exploit an analysis of the combinatorics on unions
of closed curves, arising as geodesics of such orbifold metrics.
\end{abstract}

\footnotetext[1]{\large {\em Geometric and Functional Analysis
(GAFA)}, to appear}

\maketitle

\tableofcontents

\section{To fill a circle: an introduction}
\label{one1}

Consider a compact manifold $N$ of dimension $n \geq 1$ with a
distance function $d$.  Here $d$ is not necessarily Riemannian.  The
notion of the Filling Volume, $\FillVol(N^n,d)$, of $N$ was introduced
in~\cite{Gr1}, where it is shown that when $n\geq 2$,
\begin{equation}
\label{fillvoldef}
\FillVol(N^n,d)=\inf_\gmetric{\vol(X^{n+1},\gmetric)}
\end{equation}
where $X$ is any fixed manifold such that $\partial X=N$.  Here one
can even take a cylinder $X=N\times [0,\infty)$.  The infimum is taken
over all complete Riemannian metrics $\gmetric$ on $X$ for which the
boundary distance function is $\geq d$, \ie the length of the shortest
path in $X$ between boundary points $p$ and $q$ is $\geq d(p,q)$.  In
the case $n=1$, the topology of the filling $X^2$ can affect the
infimum, as is shown by example in \cite[Counterexamples~2.2.B]{Gr1}.

The precise value of the filling volume is not known for any
Riemannian metric.  However, the following values for the canonical
metrics on the spheres (of sectional curvature $+1$) is conjectured in
\cite{Gr1}, immediately after Proposition 2.2.A.
\begin{conjecture}
\label{fillsphere}
$\FillVol(S^n,can)=\frac 1 2 \omega_{n+1}$, where $\omega_{n+1}$
represents the volume of the unit $(n+1)$-sphere.
\end{conjecture} 
This conjecture is still open in all dimensions.  The case $n=1$ can
be broken into separate problems depending on the genus of the
filling.  The filling area of the circle of length $2\pi$ with respect
to simply connected fillings (\ie by a disk) is indeed $2{\pi}$, by
Pu's inequality~\eqref{21p}, \cf \cite{Be1}, applied to the real
projective plane obtained by identifying all pairs of opposite points
of the circle.  Here one may need to smooth the resulting metric, with
an arbitrarily small change in the two invariants involved.

In Corollary~\ref{65} below, we prove the corresponding result when
the filling is by a surface of genus one.  

\begin{remark}
\label{12b}
Consider a filling by an orientable surface~$X_1$ of positive genus.
Thus $\partial X_1 = S^1$.  We apply the same technique of glueing
antipodal points of the boundary together to form a nonorientable
surface~$X_2$ (whose metric might again need to be smoothed as
before).  We are thus reduced to proving a conjectured relative
version of Pu's inequality, see discussion around formula~\eqref{16}.
This conjecture is most easily stated by passing to the orientable
double cover $X_3$ of $X_2$ with deck transformation~$\tau$, as
follows.
\end{remark}

\begin{conjecture}
\label{relPu}
Let $X$ be an orientable surface of even genus with a Riemannian
metric $\gmetric$ that admits a fixed point free, orientation
reversing, isometric involution~$\tau$.  Then there is a point $p\in
X$ with
\[
\frac {\dist _\gmetric (p, \tau(p))^2} {\area(\gmetric)} \leq
\frac{\pi}{4}.
\]
\end{conjecture} 

\begin{remark}
The analogous conjecture is false in odd genus \cite{Pa}.
\end{remark}

The original proof of Pu's inequality passed via conformal techniques
(namely, uniformisation), combined with integral geometry (see
\cite{Iv} for another proof).  It seems reasonable also to try to
apply conformal techniques to the proof of the conjectured relative
version of Pu's inequality.  This works well in the hyperelliptic
case.

Let $X$ be an orientable closed aspherical surface.  (By a surface in
this context we mean one that comes with a fixed conformal structure
and maps are conformal maps.)  

\begin{definition}
\label{63}
A {\em hyperelliptic involution\/} is a holomorphic map $J: X \to X$
satisfying $J^2=1$, with $2\genus +2$ fixed points, where $\genus$ is
the genus of $X$.  A surface $X$ admitting such an involution will be
called hyperelliptic.
\end{definition}

The involution $J$ can be identified with the nontrivial element in
the center of the (finite) automorphism group of $X$ (\cf
\cite[p.~108]{FK}) when it exists, and then such a $J$ is unique, \cf
\cite[p.204]{Mi}.  Note that the quotient map
\begin{equation}
\label{52}
Q:X\to S^2
\end{equation}
of such an involution $J$ is a conformal branched 2-fold covering.

\begin{definition}
The $2\genus+2$ fixed points of $J$ are called {\em Weierstrass
points}.  Their images in $S^2$ under the ramified double cover $Q$ of
formula~\eqref{52} will be referred to as {\em ramification points}.
\end{definition}

Our result about hyperelliptic surfaces is the following.  Given a
circle~$C_0\subset S^2$ and a point $w_0\in S^2\setminus C_0$, we can
consider the orbifold metric with equator $C_0$ and poles at $w_0$ and
at the image of $w_0$ under inversion in~$C_0$, \cf
Section~\ref{eight}.  Given a hyperelliptic surface $Q: X\to S^2$
together with a circle $C\subset X$ double covering a circle
$C_0\subset S^2$, and a point $w\in X \setminus C$, denote by
\begin{equation}
\label{63b}
AF(C, w)
\end{equation}
the pullback to $X$ of the orbifold metric for $C_0= Q(C)$ and
$w_0=Q(w)$.  Given an involution $\iota: X\to X$ fixing a circle,
denote by $X^\iota \subset X$ its fixed circle.

A closed oriented Riemann surface $X$ is called {\em ovalless real\/}
if it admits a fixed point free, antiholomorphic involution $\tau$,
see Appendix~\ref{realsurfaces} for a more detailed discussion.

\begin{theorem}
\label{relPuthm}
Let $(X,\tau,J)$ be an ovalless real hyperelliptic surface of even
genus $\genus$, and $\gmetric$ a Riemannian metric in its conformal
class.  Then there is a point $p\in X$ with
\begin{equation}
\label{55}
\frac {\dist _\gmetric (p, \tau(p))^2} {\area(\gmetric)} \leq
\frac{\pi}{4}.
\end{equation}
Specifically, there exists a curve joining $p$ and $\tau(p)$, of
length at most $\left( \frac{\pi}{4}\area(\gmetric) \right)^{1/2}$,
which consists of arcs of at most $\genus+1$ special
curves~$\gamma_i$.  Each of the $\gamma_i$ is a geodesic of the
singular constant (positive) curvature metric
\[
AF(X^{\tau \circ J}, w_i),
\]
where $w_i$ is one of the Weierstrass points, while $X^{\tau \circ J}$
is the fixed circle of the involution $\tau \circ J$.
\end{theorem}

Since every genus 2 surface is hyperelliptic \cite
[Proposition~III.7.2, page~100] {FK}, we obtain the following
corollary.

\begin{corollary}
\label{65}
All orientable genus 1 fillings of the circle satisfy the filling area
conjecture.
\end{corollary}

For a precise calculation of a related invariant called the {\em
filling radius\/} see~\cite{Ka0, KL}.  An optimal inequality for
CAT(0) metrics on the genus 2 surface is proved in \cite{KS3}.
Recently it was shown \cite{KS1} that all hyperelliptic surfaces
satisfy the Loewner inequality $\pisys_1^2 \leq \gamma_2 \area$
(see~\cite{Be2}), as well as all surfaces of genus at least 20
\cite{KS2}.  Analogues and higher dimensional optimal generalisations
of the Loewner inequality are studied in \cite{Am, BK2, IK, BCIK2},
\cf Section~\ref{two}.  Near-optimal asymptotic bounds are studied in
\cite{Ka3, BB, KS2, Sa2}.  A general framework for the study of
systolic inequalities in a topological context is proposed in
\cite{KR}.

The notion of an ovalless real hyperelliptic surface is reviewed in
Appendix~\ref{realsurfaces}.  The relation of our theorem to Pu's
inequality is discussed in Section~\ref{two}.  The two-step proof of
Theorem~\ref{relPuthm} is sketched in Section~\ref{sketch}.  An
orbifold metric on the sphere that plays a key role in the proof is
described in Section~\ref{eight}.  The key integral-geometric
ingredient of the proof appears in Section~\ref{five}.  The two steps
appear, respectively, in Section~\ref{step1} and Section~\ref{step2}.

\section{Relative Pu's way}
\label{two}
In the companion paper \cite{BCIK2}, we study certain optimal systolic
inequalities for a Riemannian manifold $(X,\gmetric)$, depending on a
pair of parameters, $n$ and~$b$.  Here $n$ is the dimension of $X$,
while~$b$ is its first Betti number.  The definitions of the systolic
invariants~$\stsys_1(\gmetric)$, $\confsys_1(\gmetric)$, and
$\pisys_1(\gmetric)$ can be found in the survey \cite{CK}, while the
Abel-Jacobi map in \cite{Li}, \cite[p.249]{Gr3}, \cf
\cite[p.~249]{Mi}, \cite[(4.3)]{BK2}, \cite{IK}.  The proof of the
inequalities involves streamlining and generalizing the techniques
pioneered in \cite{BI94, BI95} and \cite[pp.~259-260]{Gr3} (\cf \cite
[inequality~(5.14)]{CK}), resulting in a construction of Abel-Jacobi
maps from~$X$ to its Jacobi torus $\T^b$, which are area-decreasing,
with respect to suitable norms on the torus.  The norms in question
are the stable norm, the conformally invariant norm, as well as
other~$L^p$ norms.

The study of the boundary case of equality in the case $n=b+2$ depends
on the filling area conjecture and/or a suitable generalisation of
Pu's inequality, as we now discuss.  The inequality of P. Pu \cite{Pu}
(\cf \cite{Iv}) can be generalized as follows: every surface
$(S,\gmetric)$ which is not a 2-sphere satisfies the inequality
\begin{equation}
\label{21p}
\pisys_1(\gmetric)^2\leq \frac{\pi}{2}\area(\gmetric), 
\end{equation}
where the boundary case of equality is attained precisely when, on the
one hand, the surface $S$ is a real projective plane, and on the
other, the metric $\gmetric$ is of constant Gaussian curvature.  See
\cite{IK} for a more detailed discussion.

Now let $S$ be a nonorientable closed surface.  Let $\phi:\pi_1(S)\to
\Z_2$ be an epimorphism, corresponding to a map $\hat\phi:S\to \R P^2$
of absolute degree~$+1$.  We define the ``1-systole relative to
$\phi$'', denoted $\sys_1(S, \phi, \gmetric)$, of a metric $\gmetric$
on $S$, by minimizing length over loops $\gamma$ which are not in the
kernel of $\phi$:
\begin{equation}
\label{cps}
\sys_1(S, \phi, \gmetric)=\min \left\{ \length (\gamma) :
\phi([\gamma]) \not= 0\in \Z_2^{\phantom{2}} \right\}
\end{equation}

\begin{question}
\label{21}
Does every nonorientable surface $(S,\gmetric)$ and map $\hat\phi:S\to
\R P^2$ of absolute degree one, satisfy the inequality
\begin{equation}
\label{23}
\sys_1(S, \phi, \gmetric)^2\leq \tfrac{\pi}{2} \area (\gmetric),
\end{equation}
which can be thought of as a relative version of Pu's inequality?
\end{question}

This question appeared in \cite [conjecture~2.7]{CK}.  Let
\begin{equation}
\label{16}
\sigma_2=\sup_{(S,\gmetric)} \frac{\sys_1(S, \phi, \gmetric) ^2}
{\area (\gmetric)},
\end{equation}
where the supremum is over all triples $(S, \gmetric, \phi)$ as above.
Thus we ask whether $\sigma_2=\frac{\pi}{2}$.  The answer is
affirmative in the class of metrics~$\gmetric$ whose underlying
conformal structure (on the associated orientable surface) is
hyperelliptic, \cf Theorem~\ref{relPuthm}.  

Consider a genus one filling $X_1$ of the circle $S^1$, and the
associated nonorientable surface~$X_2$, \cf Remark~\ref{12b}.  Then
$X_2$ is diffeomorphic to the connected sum $3\R P^2$ of 3 copies of
the real projective plane.  A~path joining a pair of opposite points
of $S^1=\partial X_1$ corresponds to a 1-sided loop in $X_2$, and thus
defines an element of $\pi_1(X_2)$ outside the subgroup $\pi_1(X_3)
\subset \pi_1(X_2)$.  We can therefore restate Corollary~\ref{65} as
follows.

\begin{theorem}
Consider the homomorphism $\phi: \pi_1(3 \R P^2) \to \Z_2$ whose
kernel is the fundamental group of the orientable double cover of $3
\R P^2$.  Then the relative systole $\sys_1(3\R P^2, \phi, \gmetric)$
satisfies relative Pu's inequality~\eqref{23}.
\end{theorem}

An analogue for higher dimensional manifolds of the relative
1-systolic ratio was studied in \cite{Bab}.

It was shown in \cite{IK} that $\sigma_2\in [\frac{\pi}{2}, 2]$.
Furthermore, a theorem was proved in \cite{IK}, which in the case of
the manifold $\R P^2 \times \T^b$ states that every metric $\gmetric$
on $\R P^2 \times \T^b $ satisfies the following inequality:
\begin{equation}
\label{12}
\stsys_1(\gmetric)^b \pisys_1(\gmetric)^2 \leq \sigma_2
\gamma_b^{\frac{b}{2}} \vol_{b+2}(\gmetric),
\end{equation}
where $\sigma_2$ is the optimal systolic ratio from \eqref{16}.

Answering Question \ref{21} in the affirmative would allow one to
characterize the boundary case of equality in equation~\eqref{12}.
Our Theorem~\ref{relPuthm} can be thought of as a partial result in
this direction.  Calculating~$\sigma_2$ depends on calculating the
filling area of the Riemannian circle.  Removing the hyperellipticity
assumption of Theorem~\ref{relPuthm} would amount to proving the
filling area conjecture \cite{Gr1}, as explained in section
\ref{one1}.

\section{Outline of proof of Theorem~\ref{relPuthm}}
\label{sketch}

We will prove Theorem \ref{relPuthm} in a way modeled on Pu's proof.
The idea is to obtain estimates for an arbitrary metric in the
conformal class of a hyperelliptic surface, by applying integral
geometry to the quotient metric (by the hyperelliptic involution) on
the 2-sphere.  We start with a metric in the conformal class of an
ovalless real hyperelliptic surface~$(X,\tau,J)$ of even
genus~$\genus$ and area $A$, and then proceed as follows.

{\bf Step 1:} In Lemma \ref{61}, we will show that we may assume that
the commuting involutions $J$ and $\tau$ are isometries of our metric.
Hence there is an induced metric $\gmetric_1=f^2(x)\gmetric_0$ of area
$A/2$ on $S^2$ (here $f(x)$ is a nonnegative square integrable
function) which pulls back (away from the Weierstrass points) to our
original metric on $X$.  The involution $\tau$ induces an orientation
reversing isometry, denoted~$\tau_0$, on $S^2$ which can be assumed to
be the map that leaves the equator fixed and exchanges the northern
and southern hemispheres.  Clearly, $\tau_0$ maps ramification points
to ramification points, and hence there are $\genus +1$ of them in
each hemisphere.

{\bf Step 2:} There is a loop $\gamma$ on $S^2$ of $\gmetric_1$-length
$\leq \sqrt{\frac \pi 4 A}$ through a point $p$ on the equator (hence
$p$ is fixed by $\tau_0$) such that a lift of $\gamma$ starting at a
preimage, $\tilde p$, of $p$ back on $X$ is not closed.

We call a loop whose $\gmetric_1$-length $\leq \sqrt{\frac \pi 4 A}$
an {\em optimally short loop}.  The loop we find lies in a hemisphere.
In fact, we find a loop with a lift that is not closed, \ie the sum of
the winding numbers about the ramification points in the hemisphere is
odd.

Note that this proves the theorem since the endpoints of such a lift
would have to be $\tilde p$ and $\tau(\tilde p)$ since the fact that
$p$ is fixed by $\tau_0$ means that $\tau$ (which has no fixed points)
must exchange the two (distinct) preimages of $p$.

\section{Near optimal surfaces and the football}
\label{eight}

It is often useful in the proof of sharp inequalities to keep in mind
how the argument applies to spaces that are optimal or near optimal.
We present such a discusson below, and also find some optimally short
loops on these optimal surfaces.

A special role here is played by the simply connected 2-dimensional
orbifold $\AF$, which we will refer to henceforth as an (American)
football, see Figure~\ref{af}, \cf \cite{Bo}.  It is the orbifold of
constant curvature with a pair of conical singularities, each with
total angle $\pi$.  We denote the metric $\gmetric_{AF}$.

Denote by $\gmetric_0$ the round metric of Gaussian curvature $+1$ on
the 2-sphere.  Here $S^2= \C \cup \{\infty\}$, with the usual
coordinates $(x,y)= (\Re(z), \Im(z)), z\in \C$.  We have
\begin{equation}
\label{131}
\gmetric_0 = \frac{4} {\left( 1+ x^2+y^2 \right)^2} \left(
dx^2+dy^2\right),
\end{equation}
\cf \cite[p.~292]{Ri}.

\begin{proposition}
\label{71}
The football admits the following two equivalent descriptions:
\begin{enumerate}
\item
The orbifold $\AF$ is obtained from a round hemisphere by folding its
boundary in two, \ie using an identification on the boundary which is
a reflection with 2 antipodal fixed points.
\item
The rotation by $\pi$ around the $z$-axis in $\R^3$ gives rise to a
ramified double cover $D: (S^2,\bar \gmetric_0) \to \AF$, which is a
local isometry away from the poles.
\end{enumerate}
Furthermore, the orbifold metric $\gmetric_{AF}$ can be expressed in
terms of the standard metric by means of the following conformal
factor:
\begin{equation}
\label{4.2}
\gmetric_{AF}= \frac{1+ {\mathcal O} (r)}{4r} \gmetric_0,
\end{equation}
where ${\mathcal O}(r) \to 0$ as $r\to 0$, while $r=| z | =
\sqrt{x^2+y^2}$.
\end{proposition}

\begin{proof}
The equivalence of the first two definitions is clear.  Let us relate
the first definition to the third.  We represent a hemisphere defining
the orbifold as the upperhalf plane $H^2\subset \C$.  Then the
orbifold can be viewed as a quotient of $(H^2, \gmetric_0)$, where the
positive and negative rays of the $x$-axis are identified.  Denote by
$w$ the complex coordinate in~$H^2$.

By the uniqueness of the underlying conformal structure, there exists
a conformal map $c: (H^2,w) \to (S^2,z)$, identifying the orbifold
with the 2-sphere, with complex coordinate $z$.  In our coordinates,
we can represent $c$ by the map $z= w^2$.  Here the singular points of
the orbifold correspond to $0$ and~$\infty$.  Since $\frac {dc} {dw} =
2w$, the metric $\gmetric_0$ pulls back under~$c$ to the quadratic
differential $c^*(\gmetric_0)(w) = \mu(|w|^2) \gmetric_0$, where
\[
\mu(t) = 4t \frac{(1+t)^2}{(1+t^2)^2} = 4t (1 + {\mathcal O}(t)).
\]
Note that $\gmetric_0$ restricted to $H^2$ is precisely the orbifold metric
$c^*(\gmetric_{AF})$ we are looking for.  Thus
\[
\begin{aligned}
c^*(\gmetric_{AF}) & = \frac{1}{\mu(|w|^2)} c^*(\gmetric_0) \cr & =
\frac{1}{\mu(|z|)} c^*(\gmetric_0) \cr & = c^*\left( \frac{1+{\mathcal
O}(|z|)}{4{|z|}} (\gmetric_0) \right),
\end{aligned}
\]
as required.  
\end{proof}

\begin{remark}
\label{82}
The function $\frac{1}{\sqrt{r}}$ is locally square integrable
in~$\C$, which is important for our applications, \cf
Proposition~\ref{73} of Appendix~\ref{bee}.
\end{remark}

\begin{example}
The football $\AF$ arises in the typical example where the associated
systolic ratio is close to being optimal.  Thus, start with a round
$\R P^2$.  Attach a small handle.  The orientable double cover $X$ can
be thought of as the unit sphere in $\R^3$, with two little handles
attached at north and south poles, \ie at the two points where the
sphere meets the $z$-axis.  Then one can think of the hyperelliptic
involution $\hyperelliptic$ as the rotation of $X$ by $\pi$ around the
$z$-axis.  The six Weierstrass points are the six points of
intersection of $X$ with the $z$-axis.  The orientation reversing
involution $\tau$ on $X$ is the restriction to~$X$ of the antipodal
map in $\R^3$.  The composition $\tau\circ \hyperelliptic$ is the
reflection fixing the $xy$-plane, in view of the following matrix
identity:
\begin{equation}
\label{92b}
\begin{pmatrix}
-1 & 0 & 0 \cr 0 & -1 & 0 \cr 0 & 0 & -1
\end{pmatrix}
\begin{pmatrix}
-1 & 0 & 0 \cr 0 & -1 & 0 \cr 0 & 0 & 1
\end{pmatrix}
=
\begin{pmatrix}
1 & 0 & 0 \cr 0 & 1 & 0 \cr 0 & 0 & -1
\end{pmatrix}
.
\end{equation}
Meanwhile, the induced orientation reversing involution $\tauzero$ on
$S^2$ can just as well be thought of as the reflection in the
$xy$-plane.  This is because, at the level of the 2-sphere, it is
``the same thing as'' the composition $\tau\circ \hyperelliptic$.
Thus the fixed circle of $\tauzero$ is precisely the equator, \cf
formula~\eqref{72a}.  Then one gets a quotient metric on $S^2$ which
is roughly that of the western hemisphere, with the boundary longitude
folded in two, \cf Proposition~\ref{71}(2).  The metric has little
bulges along the $z$-axis at north and south poles, which are
leftovers of the small handle.
\end{example}

\section{Finding a short figure eight geodesic}
\label{five}

\begin{figure}
\includegraphics[height=9cm]{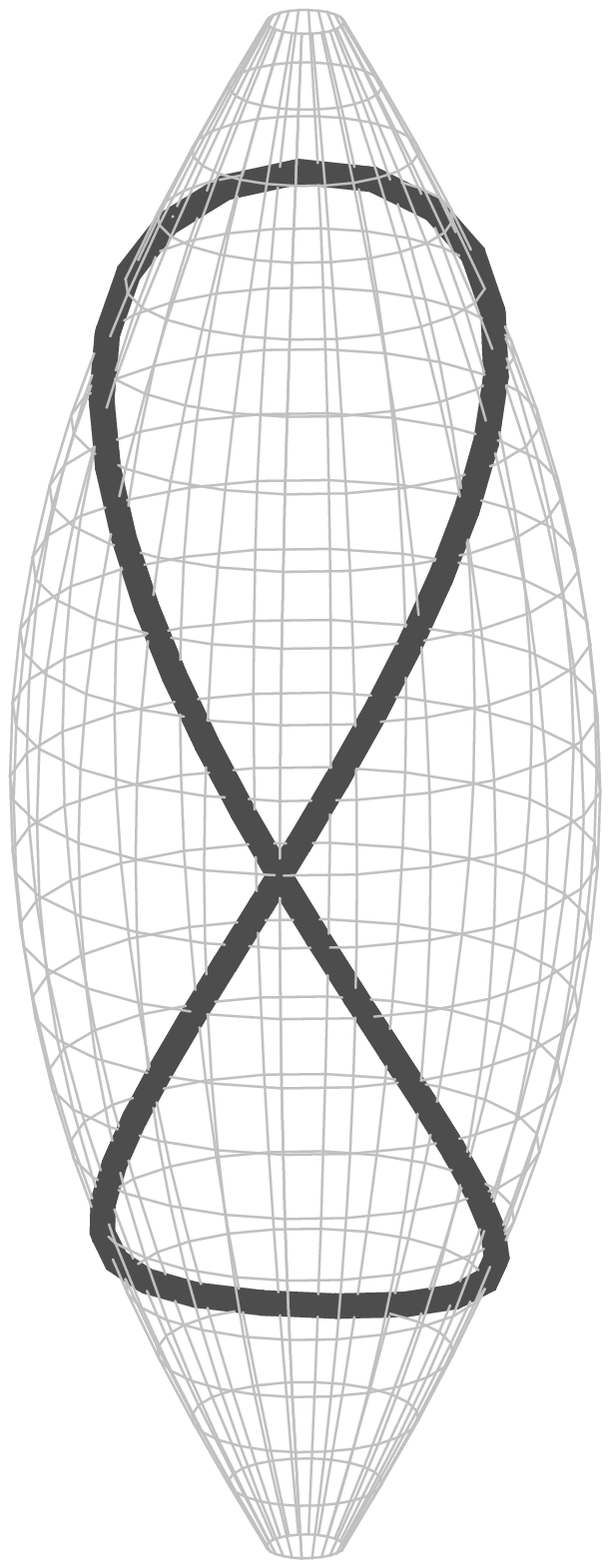}
\includegraphics[height=9cm]{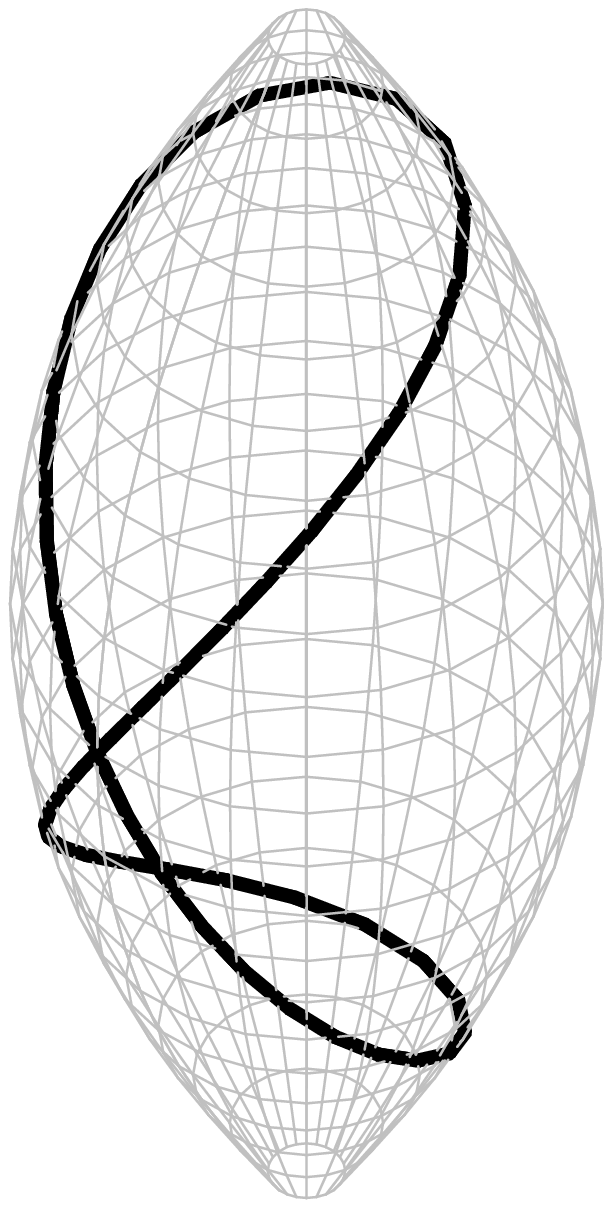}
\includegraphics[height=9cm]{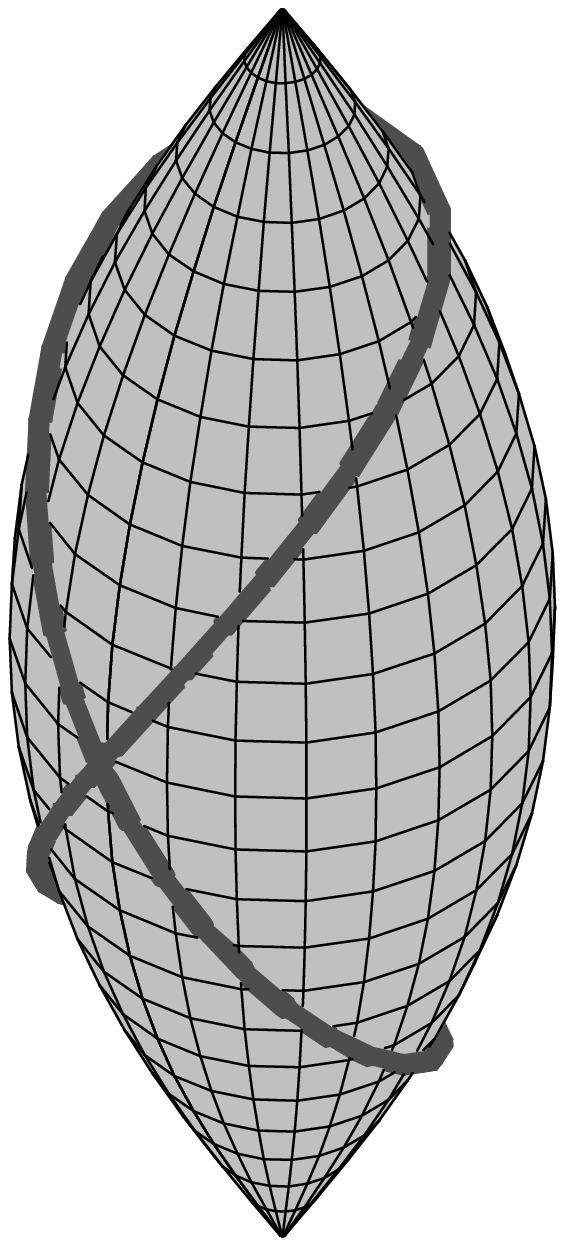}
\caption{Figure eight geodesics on a football, \cf \eqref{4.2}}
\label{af}
\end{figure}

In our case, the optimal metric is not round, so it is natural to
consider geodesics on the football $\AF$ instead of the round sphere.
The branched double cover $D:S^2\to \AF$ of Proposition~\ref{71}(2)
pulls back $\gmetric_{AF}$ to the standard metric $\gmetric_0$ on
$S^2$.  This allows us easily to see the shape of the
$\gmetric_{AF}$-geodesics.  They are of figure eight type with their
point of self intersection on the equator (except the ones through the
poles), see Figure~\ref{af}.  Note also that the loop in each of the
northern and southern hemispheres winds around the pole precisely
once.

Starting with a square integrable function $f(\sigma)$ on $\AF$, we
can pull it back to $S^2$ and apply Lemma~\ref{73} to obtain the
following.

\begin{lemma}
\label{intgeomfootball}
There is a $\gmetric_{AF}$-geodesic $\gamma$ satisfying the following
estimate:
\[
\left( \int_\gamma f(\gamma(t))dt \right) ^2\leq \pi
\int_{S^2}f^2(\sigma)d\sigma = 2\pi\area(f^2\gmetric_{AF}).
\]
\end{lemma}

The factor of 2 comes from the fact that the area of the pull back of
the (possibly singular) metric $f^2\gmetric_{AF}$ in $S^2$ (under the
double cover~$S^2 \to \AF$) has twice the area of the original metric.
Lemma~\ref {intgeomfootball} implies the following.  Assume that our
function $f$ is symmetric with respect to the inversion in the
equator.  Then with respect to the metric~$f^2\gmetric_{AF}$, the
length $L$ of each of the two hoops of the figure eight curve,
satisfies the desired inequality $L\leq \sqrt{\frac {\pi} {2}
A(f^2\gmetric_{AF})}.$ This inequality is sharp, since
$A(f^2\gmetric_{AF})$ is half the area $A$ of the original
hyperelliptic surface.

\section{Proof of circle filling: Step 1}
\label{step1}

Let $X$ be an ovalless real hyperelliptic surface of even genus
$\genus$.  Let $\hyperelliptic: X\to X$ be the hyperelliptic
involution, \cf Definition~\ref{63}.  Let $\tau$ be the
antiholomorphic involution defining the real structure, \cf
Appendix~\ref{realsurfaces}.

\begin{lemma}
\label{92}
The involution $\tau$ commutes with $J$ and descends to $S^2$.  The
induced involution $\tauzero: S^2\to S^2$ is an inversion in a circle
$C_0= Q(X^{\tau\circ J})$.  The set of ramification points is
invariant under the action of $\tauzero$ on $S^2$.
\end{lemma}

\begin{proof}
By the uniqueness of $J$ (\cf Appendix~\ref{realsurfaces}), we have
the commutation relation
\begin{equation}
\label{commutation}
\tau \circ \hyperelliptic = \hyperelliptic \circ \tau,
\end{equation}
\cf relation~\eqref{92b}.  Therefore $\tau$ descends to an involution
$\tauzero$ on the sphere.  Suppose~$\tauzero$ were fixed point free.
Then it would be conjugate to the antipodal map of~$S^2$.  Therefore,
the set of the $2\genus+2$ ramification points on $S^2$ is centrally
symmetric.  Since there is an odd number, $\genus+1$, of ramification
points in a hemisphere, a generic great circle $A\subset S^2$ has the
property that its inverse image $Q^{-1}( A ) \subset X$ is connected.
Thus both involutions $\tau$ and $\hyperelliptic$, as well as
$\tau\circ \hyperelliptic$, act fixed point freely on the loop
$Q^{-1}( A) \subset X$, which is impossible.  Therefore $\tauzero$
must fix a point.  It follows that $\tauzero$ is an inversion in a
circle.
\end{proof}

The next lemma allows us to assume that $J$ and $\tau$ are isometries
of the metric.

\begin{lemma}
\label{61}

Let $F$ denote a conformal involution of $X$, and let $\bar \gmetric
=\frac 1 2 \{ \gmetric +F^* \gmetric \}$ be the averaged metric. Then
for all $p,q \in X$ we have
\begin{equation} 
\label{81}
\frac {\dist _{\bar \gmetric} (p, q )^2} {\area(\bar \gmetric)}
\geq\frac{ \left( \frac{ 1}{ 2}\{\dist _\gmetric (p, q)+\dist
_\gmetric (F(p), F(q))\} \right) ^2} {\area(\gmetric)}.
\end{equation} 
and hence if Theorem \ref{relPuthm} holds for $\bar \gmetric$ then it
holds for $\gmetric$, as well.
\end{lemma}
\begin{remark}
Applying Lemma \ref{61} to $F=\tau$ and $q=\tau(p)$, we obtain
a~$\tau$-invariant metric $ \gmetric^-$, conformal to $\gmetric$, such
that
\[
\frac {\dist _{ \gmetric ^-} (p, q )^2} {\area( \gmetric ^-)} \geq
\frac {\dist _{ \gmetric} (p, q )^2} {\area( \gmetric)}.
\]

If we apply Lemma \ref{61} to $\gmetric^-$ and $F=J$, $q=\tau(p)$, we
obtain a~$J$-invariant metric $\bar \gmetric$, also conformal
to~$\gmetric$. Now inequality~\eqref{81} implies that if Theorem
\ref{relPuthm} holds for $\bar \gmetric$, then it holds for $\gmetric$
as well.  Moreover, since $J$ and $\tau$ commute, $ \bar \gmetric$ is
also $\tau$-invariant.
\end{remark}

\begin{proof}
Since $F$ is a conformal diffeomorphism, we have $F^* \gmetric =f^2
\gmetric$ for some smooth positive $f$.  Then $\bar \gmetric =\frac 1
2 (f^2+1)\gmetric$.  (Note that a similar proof works also with $\bar
\gmetric = \left( \frac 1 2 (f+1) \right) ^2 \gmetric$ as an average.)
Let $dx$ represent the volume element of $\gmetric$.  We have
\[
\begin{aligned}
\area(\bar \gmetric) & =\int \frac 1 2 (f^2(x)+1)dx \cr & =\frac 1 2
\left( \int f^2(x)dx+\int 1 dx \right) \cr &= \frac 1 2 (\area
(F^*(\gmetric))+\area(\gmetric))\cr &=\area(\gmetric),
\end{aligned}
\]
and thus the denominator is unchanged.  As for the numerator, start
with any curve $\gamma$ from $p$ to $q$, then $F(\gamma)$ is a curve
from $F(p)$ to the point $F(q)$.  The length element for $\bar
\gmetric$ is $ \sqrt{\frac{f^2 + 1}{2}} \geq \frac{f + 1}{2} ,$ so
that length of the curve $\gamma$ for the averaged metric is bounded
below by the average of the length of the curves $\gamma$ and $F \circ
\gamma$ for the original metric.  Thus taking $\gamma$ to be a $\bar
\gmetric$ minimizing path, we obtain
\[
\begin{aligned}
\dist _{\bar \gmetric} (p, q) & = L_{\bar \gmetric}(\gamma) \cr & \geq
\frac 1 2 (L_{\gmetric}(\gamma)+L_{\gmetric}(F \circ \gamma))\cr &
\geq \frac 1 2\{\dist _\gmetric (p, q)+\dist _\gmetric (F(p),
(F(q)))\},
\end{aligned}
\]
proving inequality \eqref{81}.
\end{proof}

Now if inequality~\eqref{55} holds for $\bar \gmetric$ then there is a
$p$ such that $\frac \pi 4 \geq\frac {\dist _{\bar \gmetric} (p,
\tau(p))^2} {\area(\bar \gmetric)}$.  Hence either $\frac \pi 4
\geq\frac {\dist _{\gmetric} (p, \tau(p))^2} {\area(\gmetric)}$, or
$\frac \pi 4 \geq\frac {\dist _{\gmetric} (F(p), \tau(F(p)))^2}
{\area(\bar \gmetric)}$, and the inequality follows for $\gmetric$, as
well.  From now on we will assume that~$J$ and~$\tau$ are isometries.

\section{Proof of  circle filling: Step 2}
\label{step2}

Let $\Ci=Q(X^{\tau\circ J}) \subset S^2$ be the fixed circle of the
involution~$\tauzero$.  Let $H\subset S^2$ be a connected component of
$S^2\setminus \Ci$.  Since the set of Weierstrass points is invariant
under $\tau$, there is an odd number of ramification points, namely
$\genus+1$, in the ``hemisphere'' $H$.

We need to find a optimally short loop (\ie a loop satisfying $L\leq
\sqrt{\frac \pi 4 A}$) based at a point on the ``equator'' $\Ci$ whose
lift does not close up, and hence connects an orbit of $\tau$.

\begin{lemma}
\label{121a}
For every ramification point $w\in H$, there is an optimally short
loop (half of a figure 8 loop) in $H$ with unit winding number around
$w$.
\end{lemma}

\begin{proof}
By using a Mobius transformation if necessary, we can assume
that~$\Ci$ is a great circle and our ramification point $w\in S^2$ is
the polar point of the great circle $\Ci$.

Recall that $S^2$ carries a metric induced from $X$.  We then can use
a conformal diffeomorphism from $\AF$ to $S^2$ which takes the
ramification points to the poles, to pull the metric back to $\AF$.

Namely, we compare our surface $X$ to $AF(X^{\tau\circ J}, w)$, \cf
formula~\eqref{63b}.  Here $AF(X^{\tau\circ J}, w)$ double covers
$\AF$, while the circle $X^{\tau\circ J}$ double covers $\Ci$ under
$Q$.  Recall that the conformal factor of our metric may be assumed
invariant under~$\tau$.  Hence Lemma \ref {intgeomfootball} produces
an optimally short loop (half of a figure 8 curve) with unit winding
number around $w$, based at a point of $\Ci$.

By placing ourselves in a generic situation, we can assume that the
loops we work with do not pass through ramification points, and thus
the winding number is well defined.  Namely, we assume that the
geodesic we pick on $\AF$ is generic.  This worsens the resulting
estimate by an arbitrarily small $\epsilon>0$.  But if we prove that
the optimal systolic ratio is within $\epsilon$ of $\pi/2$ for
arbitrary $\epsilon$, then it certainly follows that it equals
$\pi/2$.
\end{proof}

Note in the above we will need to take a different Mobius
transformation for each $w$.

If a loop produced by Lemma~\ref{121a} encircles no other ramification
points, then its lift to the ovalless real hyperelliptic Riemann
surface~$(X,\tau,J)$ connects an orbit of $\tau$, and we are done with
Step 2.  Otherwise we argue as follows.

Fix a finite set $S$ of points in the hemisphere.  In the end, $S$
will be the set of ramification points, but the argument involves
subsets, too.  The loops involved in the argument are not assumed
simple.  Given a loop $b$, denote by $W_S(b)$ the set of points of $S$
for which the winding number of $b$ is odd.

\begin{lemma}
\label{171}
Let $b$ and $c$ be smooth optimally short loops in a disk, such that
the set $|W_S(c)\setminus W_S(b)|$ is odd.  Then there is an optimally
short loop $h$ (with the same basepoint as either $b$ or $c$) such
that $|W_S(h)|$ is odd.
\end{lemma}

\begin{proof}
Consider a double ramified cover of the disk, with ramification points
at~$S$.  Given a loop, the number of points with odd winding number is
even if and only if the lift of the loop to the cover closes up.

Choose a fixed lift of the loop $b$, and a fixed lift of the loop $c$.
Suppose both lifts are closed.  The points of the intersection $b\cap
c$ fall into two types: ``crossroads'', where their lifts intersect,
as well; and ``bridges'', where the lifts do not intersect, \ie lie in
different sheets of the cover.

It suffices to find two intersection points of different types.
Indeed, if $p$ and $q$ are such points, consider the arc of $b$
between $p$ and $q$ not containing the basepoint, and the similar arc
of $c$. Exchanging these arcs produces a new pair of loops of the same
combined length, whose lifts do not close up.  Hence one of them is
optimally short, as required.

We need the fact that a self-intersection of one of the loops remains
a self-intersection after lifting.  Suppose it does not.  Then we have
a sub-loop, not containing the basepoint, whose lift is not closed.
Remove this sub-loop from the original loop.  Then the remaining loop
has the desired property: it is optimally short, and its lift is not
closed.

To find a pair of intersection points of different types, we argue by
contradiction.  Assume that they are all crossroads (if bridges,
change the lift of~$b$).  Thus, we have a (homeomorphic) lift of the
union of the loops $b\cup c$.  Hence every simple loop contained in
this union lifts to a closed curve, and therefore encircles an even
number of ramification points.  Consider the graph $\Gamma$ in the
disk, defined by the union of the two loops.  Let~$F_i$ be the faces
(\ie connected components of the complement) of~$\Gamma$ containing
points from $W_S(c) \setminus W_S(b)$.  Consider the subgraph~$\cup_i
F_i$ of~$\Gamma$, and note that
\begin{equation}
\label{61b}
W_S(c)\setminus W_S(b) \subset \cup_i F_i.
\end{equation}
Since $|W_S(c)\setminus W_S(b)|$ is odd, formula~\eqref{61b} implies
that one of the faces, say $F_{i_0}$, must contain an odd number of
points $w_j$ from $W_S(c)\setminus W_S(b)$.  By construction, for each
such $w_j\in F_{i_0}$, the loop $c$ has odd winding number, while the
loop~$b$ has even winding number.  Any ramification point $p\in
F_{i_0}$ can be connected to a $w_j$ by a path disjoint from $b\cup
c$.  Hence the loops $b$ and $c$ have the same winding number with
respect to~$p$ as with respect to $w_j$.  Hence in fact we have $p\in
W_S(c)\setminus W_S(b)\subset S$.

We conclude that $F_{i_0}$ contains an odd number of ramification
points.  Hence the lift of its boundary does not close up.  This
contradicts the assumption that all intersection points are
crossroads, and proves the lemma.
\end{proof}

\begin{lemma}
\label{ind}
Let $S$ be an odd pointset in the interior of a Riemannian disk.
Consider a family of smooth simple loops $\gamma_1, \ldots, \gamma_k$
based at distinct points of the boundary.  Let~$L>0$, and assume each
loop is of length at most $L$, and that each point of~$S$ is in the
interior of at least one of the loops.  Then the union of the loops
$\gamma_i$ contains a subloop $\gamma$, also based at a point of the
boundary and of length at most $L$, such that~$W_S(\gamma)$ is odd.
\end{lemma}

\begin{proof}
Let $|S|=2n+1$.  The proof is by induction on $n$.  Let~$b=\gamma_1$.
If~$W_S(b)$ is odd, the desired short loop is $\gamma=b$.
Assume~$W_S(b)$ is even.  By the inductive hypothesis applied to the
odd set
\begin{equation}
S_0= S \setminus W_S(b)
\end{equation}
in place of $S$, there is a short loop $c$ such that $W_{S_0}(c)$ is
odd.  {\em A priori\/} the loop $c$ might have odd winding number with
respect to additional points in $W_S(b)$, but at any rate the equality
\[
W_S(c) \setminus W_S(b) = W_{S_0}(c)
\]
allows us to apply Lemma~\ref{171} to the pair of loops $b$ and $c$,
completing the proof of Lemma~\ref{ind} and Step 2.
\end{proof}

\appendix

\section{Ovalless reality and hyperellipticity}
\label{realsurfaces}

\subsection{Hyperelliptic surfaces}
For a treatment of hyperelliptic surfaces, see \cite[p.~60-61]{Mi}.
By \cite[Proposition 4.11, p. 92]{Mi}, the affine part of a
hyperelliptic surface~$X$ is given by a suitable equation of the form
\begin{equation}
\label{a1}
w^2 = f(z)
\end{equation}
in $\C^2$, where $f$ is a polynomial.  On this affine part, the map
$J$ is given by $J(z,w) = (z,-w)$, while the hyperelliptic quotient
map $Q:X \to S^2$ is represented by the projection onto the
$z$-coordinate in $\C^2$.

A slight technical problem here is that the map $X\to \C P^2$, whose
image is the compactification of the solution set of \eqref{a1}, is
not an imbedding.  Indeed, there is only one point at infinity, given
in homogeneous coordinates by $[0:w:0]$.  This point is a singularity.
A way of desingularizing it using an explicit change of coordinates at
infinity is presented in \cite[p.~60-61]{Mi}.  The resulting smooth
surface is unique~\cite[Theorem, p.~100]{DS}.

To explain what happens ``geometrically'', note that there are two
points on our affine surface ``above infinity''.  This means that for
a large circle $|z|=r$, there are two circles above it satisfying
equation~\eqref{a1} where $f$ has even degree $2\genus+2$ (for a
Weierstrass point we would only have one circle).  To see this,
consider $z=re^{ia}$.  As the argument $a$ varies from 0 to $2\pi$,
the argument of $f(z)$ will change by $(2\genus+2)2\pi$.  Thus, if
$(re^{ia},w(a))$ represents a continuous curve on our surface, then
the argument of $w$ changes by $(2\genus+2)\pi$, and hence we end up where
we started, and not at $-w$ (as would be the case were the polynomial
of odd degree).  Thus there are {\em two\/} circles on the surface
over the circle~$|z|=r$.  We conclude that to obtain a smooth compact
surface, we will need to add two points at infinity, \cf discussion
around \cite[formula~(7.4.1), p.~102]{FK}.

Thus, the affine part of $X$, defined by equation~\eqref{a1}, is a
Riemann surface with a pair of punctures $p_1$ and $p_2$.  A
neighborhood of each~$p_i$ is conformally equivalent to a punctured
disk.  By replacing each punctured disk by a full one, we obtain the
desired compact Riemann surface~$X$.  The point at infinity
$[0:w:0]\in \C P^2$ is the image of both~$p_i$.

\subsection{Ovalless surfaces}
In our case, the hyperelliptic Riemann surface $X$ admits an
antiholomorphic involution $\tau$.  In the literature, the components
of the fixed point set of~$\tau$ are sometimes referred to as
``ovals''.  Since in our situation $\tau$ is fixed point free, we
introduce the following terminology.
\begin{definition}
\label{real}
A hyperelliptic surface $(X,J)$ of even genus $\genus>0$ is called
{\em ovalless real\/} if one of the following equivalent conditions is
satisfied:
\begin{enumerate}
\item
$X$ admits an imaginary reflection, \ie a fixed point free,
orientation reversing involution $\tau$;
\item
the affine part of $X$ is the locus in $\C^2$ of the equation
\begin{equation}
\label{71a}
w^2=-P(z),
\end{equation}
where $P$ is a monic polynomial, of degree~$2\genus+2$, with real
coefficients, no real roots, and with distinct roots.
\end{enumerate}
\end{definition}

\begin{lemma}
\label{realhyper}
The two ovalless reality conditions of Definition~\ref{real} are
equivalent.
\end{lemma}

\begin{proof}
A related result appears in \cite[p.~170, Proposition~6.1(2)]{GrH}.
To prove the implication $2\Longrightarrow 1$, note that complex
conjugation leaves the equation invariant, and therefore it also
leaves invariant the locus of \eqref{71a}.  A fixed point must be
real, but $P$ is positive hence~\eqref{71a} has no real solutions.
There is no real solution at infinity, either, as there are two points
at infinity, which are not Weierstrass points since~$P$ is of even
degree, as discussed above.  The desired imaginary reflection~$\tau$
switches the two points at infinity, and, on the affine part of the
Riemann surface, coincides with complex conjugation $(z,w)\mapsto
(\bar z, \bar w)$ in~$\C^2$.

To prove the implication $1\Longrightarrow 2$, note that by
Lemma~\ref{92}, the induced involution $\tauzero$ on $S^2= X/
\hyperelliptic$ may be thought of as complex conjugation, by choosing
the fixed circle of $\tauzero$ to be the circle
\begin{equation}
\label{72a}
\R\cup \{\infty\} \subset \C\cup \{\infty\} = S^2.
\end{equation}

Since the surface is hyperelliptic, it is the smooth completion of the
locus in~$\C^2$ of {\em some\/} equation of the form \eqref{71a}, as
discussed above.  Here $P$ is of degree~$2\genus+2$ with distinct
roots, but otherwise to be determined.  The set of roots of $P$ is the
set of (the $z$-coordinates of) the Weierstrass points.  Hence the set
of roots must be invariant under $\tauzero$.  Thus the roots of the
polynomial either come in conjugate pairs, or else are real.
Therefore $P$ has real coefficients.  Furthermore, the leading
cofficient of $P$ may be absorbed into the $w$-coordinate by
extracting a square root.  Here we may have to rotate $w$ by $i$, but
at any rate the coefficients of $P$ remain real, and thus $P$ can be
assumed monic.

If $P$ had a real root, there would be a ramification point fixed by
$\tauzero$.  But then the corresponding Weierstrass point must be
fixed by $\tau$, as well!  This contradicts the fixed point freeness
of $\tau$.  Thus all roots of~$P$ must come in conjugate pairs.
\end{proof}

\section{A double fibration and integral geometry}
\label{bee}

In the proof of Pu's theorem, as well as our argument here (see
Lemma~\ref{intgeomfootball}), one uses the following integral geometry
result.  It has its origin in results of P. Funk \cite{Fu} determining
a symmetric function on the two-sphere from its great circle
integrals, \cf \cite[Proposition~2.2, p.~59]{He}, as well as Preface
therein.

\subsection{A double fibration of $SO(3)$}

Note that the sphere $S^2$ is the homogeneous space of the Lie group
$SO(3)$, so that $S^2= SO(3)/SO(2)$.  Denote by $SO(2)_\sigma$ the
fiber over (stabilizer of) a typical point $\sigma\in S^2$.  The
projection
\begin{equation}
\label{pp}
p: SO(3)\to S^2
\end{equation}
is a Riemannian submersion for the standard metric of constant
sectional curvature $\frac{1}{4}$ on $SO(3)$.  The total space $SO(3)$
admits another Riemannian submersion, which we denote
\begin{equation}
\label{qq}
q: SO(3)\to \widetilde{\R P^2},
\end{equation}
whose typical fiber~$\gamma$ is an orbit of the geodesic flow on
$SO(3)$ viewed as the unit tangent bundle of $S^2$.  

\renewcommand{\arraystretch}{1.3}
\begin{figure}
\[
\xymatrix { SO(1)\ar[dr]&& \left(\widetilde{ \R P^2}, d\gamma \right)
\\ & SO(3)\ar[ur]^q \ar[dr]^p& \\ \gamma\quad\ar[ur] \ar[rr]_{
\hbox{great\ circle} } & & (S^2, d\sigma) }
\]
\caption{\textsf{Integral geometry on $S^2$, \cf \eqref{pp} and
\eqref{qq}}}
\label{fig1}
\end{figure}
\renewcommand{\arraystretch}{1}

Each orbit $\gamma$ projects under $p$ to a great circle on the
sphere.  We think of the base space $\widetilde{\R P^2}$ of $q$ as the
configuration space of oriented great circles on the sphere, with
measure $d\gamma$.  The diagram of Figure~\ref{fig1} illustrates the
maps defined so far.

\subsection{Integral geometry on $S^2$} 

Let $f$ be a square integrable function on $S^2$, which is positive
and continuous except possibly for a finite number of points where $f$
either vanishes or has a singularity of type~$\frac{1}{\sqrt{r}}$, \cf
Remark~\ref{82}.

\begin{proposition}
\label{73}
There is a great circle $\gamma$ such that the following two
equivalent inequalities are satisfied:

\begin{enumerate}
\item
we have $\left(\int_\gamma f(\gamma(t))dt \right)^2\leq \pi
\int_{S^2}f^2(\sigma)d\sigma$, where $t$ is the arclength parameter
and $d\sigma$ is the Riemannian measure on the standard unit sphere.
\item

there is a great circle of length $L$ in the metric $f^2 \gmetric_0$, so that
$L^2\leq \pi A$, where $A$ is the Riemannian surface area of the
metric $f^2 \gmetric_0$.
\end{enumerate}
In the boundary case of equality in either inequality, the function
$f$ must be constant.
\end{proposition}

\begin{proof}
The proof is an averaging argument and shows that the average length
of great circles is short.

Denote by $E_\gmetric(\gamma)$ the energy, and by $L_\gmetric(\gamma)$
the length, of a curve~$\gamma$ with respect to a (possibly singular)
metric $ \gmetric= f^2 \gmetric_0.$ We apply Fubini's theorem
\cite[p.~164]{Ru} twice, to both $p$ and $q$, to obtain
\[
\begin{aligned}
\frac{\left(\int _{\widetilde{\R P^2}} L_\gmetric(\gamma) ^2\right)
}{2 \pi} & \leq \int _{\widetilde{\R P^2}} E_\gmetric(\gamma) d\gamma
\cr & = \int_{\widetilde{ \R P^2}} d\gamma \left( \int_\gamma f^2\circ
p \circ \gamma(t) dt \right) \cr & = \int_{SO(3)} f^2 \circ p \cr & =
\int_{S^2} f^2 \; d\sigma \left( \int _{SO(2)_\sigma} 1 \right) \cr &
= 2\pi \; \area(\gmetric),
\end{aligned}
\]
proving the formula since $\area(\widetilde{\R P^2})=4\pi$.  In the
boundary case of equality, length and energy must be equal and
therefore the function $f$ must be constant along every great circle,
hence constant everywhere on $S^2$.
\end{proof}

\section*{Acknowledgments}
The authors are grateful to Hershel Farkas, Shmuel Krushkal, and Steve
Shnider for helpful discussions of hyperelliptic surfaces; to Mathieu
Dutour for a helpful discussion of the combinatorial arguments of
section~\ref{step2}; and to Franz Auer for the Super Bowl figures of
section \ref{five}.


\end{document}